\input amssym.def
\input amssym
\magnification=1200
\parindent0pt
\hsize=16 true cm \baselineskip=13  pt plus .2pt $ $

\def\Z{\Bbb Z}
\def\A{\Bbb A}
\def\S{\Bbb S}
\def\D{\Bbb D}
\def\R{\Bbb R}

\centerline {\bf  On finite groups acting on homology 4-spheres}

\centerline {\bf  and finite subgroups of ${\rm SO}(5)$}

\bigskip \bigskip

\centerline {Mattia Mecchia and  Bruno Zimmermann}

\bigskip

\centerline {Universit\`a degli Studi di Trieste} \centerline {Dipartimento di
Matematica e Informatica} \centerline {34127 Trieste, Italy}

\bigskip \bigskip

Abstract.  {\sl  We show that a finite group which admits a faithful, smooth,
orientation-preserving action on a homology 4-sphere, and in particular on the
4-sphere, is isomorphic to a subgroup of the orthogonal group ${\rm SO}(5)$, by
explicitly determining the various groups which can occur (up to an
indetermination of index two in the case of solvable groups). As a consequence
we obtain also a characterization of the finite groups which are isomorphic to
subgroups of the orthogonal groups ${\rm SO}(5)$ and ${\rm O}(5)$.}

\bigskip \bigskip

{\bf 1. Introduction}

\medskip

We are interested in the class of finite groups which admit an
orientation-preserving action on a homology 4-sphere, and in particular on the
4-sphere $S^4$, and to compare it with the class of finite subgroups of the
orthogonal group ${\rm SO}(5)$ acting linearly on $S^4$. As a general
hypothesis, all actions in the present paper will be faithful,
orientation-preserving and smooth (or locally linear).

\medskip

By the recent geometrization of finite group actions on 3-manifolds, every
finite group action on the 3-sphere is conjugate to an orthogonal action; in
particular, the finite groups which occur are exactly the well-known finite
subgroups of the orthogonal groups ${\rm SO}(4)$ or ${\rm O}(4)$. Finite groups
acting on arbitrary homology 3-spheres are considered in [MeZ2] and [Z]; here
some other finite groups occur and the situation is still not completely
understood, see section 4  for  a short discussion of the situation in
dimension three.

\medskip

It is no longer true that finite group actions on the 4-sphere are conjugate to
orthogonal actions; for example, it is well-known that the Smith conjecture
does not hold in dimension four (that is, the fixed point set of a periodic
diffeomorphism of the 4-sphere may be a knotted 2-sphere). However, the class
of finite groups admitting an action on the 4-sphere, and more generally also
on a homology 4-sphere, should still coincide with the class of finite
subgroups of the orthogonal group ${\rm SO}(4)$. Up to an indetermination of
index two in the case of solvable groups, this is a consequence of the
following main result of the present paper which explicitly determines the
various groups which can occur.

\bigskip

{\bf  Theorem.} {\sl  A finite group $G$ which admits an orientation-preserving
action on a homology 4-sphere is isomorphic to one of the following groups:

\medskip

i) a subgroup of the Weyl group $W = (\Z_2)^4 \rtimes \S_5$;

\medskip

ii) $\A_5$, $\S_5$, $\A_6$ or $\S_6$;

\medskip

iii)  an orientation-preserving subgroup of ${\rm O}(3) \times {\rm O}(2)$;

\medskip

iv)  a subgroup of ${\rm SO}(4)$, or a 2-fold extension of such a group.

\medskip

If $G$ is nonsolvable (but most likely also in general), iv) can be replaced by
the stronger:

\smallskip

iv')  a subgroup of ${\rm O}(4)$.}

\bigskip

Note that the different cases of the Theorem are not mutually exclusive. Since all
groups in the Theorem (except maybe the 2-fold extensions in iv)) are subgroups of
${\rm SO}(5)$ we have the following (cf. [E], Question 6 and Problem 9):

\bigskip

{\bf  Corollary 1.} {\sl  A finite group $G$ which admits an
orientation-preserving action on a homology 4-sphere is isomorphic to a
subgroup of ${\rm SO}(5)$, or possibly, if $G$ is solvable, to a 2-fold
extension of a subgroup of ${\rm SO}(4)$. }

\bigskip

The group $W$ in the Theorem can be described as follows. Consider the semidirect
product $\tilde W = (\Z_2)^5 \rtimes \S_5$ where the symmetric group $\S_5$
acts on the normal subgroup $(\Z_2)^5$ by permuting the components (Weyl group
or wreath product $\Z_2 \wr \S_5$). This semidirect product acts orthogonally
on euclidean 5-space by inversion and permutation of coordinates, and its
subgroup of index two of orientation-preserving elements is the semidirect
product $W = (\Z_2)^4 \rtimes \S_5$ in the Theorem (composing the
orientation-reversing elements in the original action of $\S_5$ on $\R^5$ with
${\rm -id}_{\R^5}$); note that $\tilde W = W \oplus \Z_2$,  with $\Z_2 =
\langle {\rm -id}_{\R^5} \rangle$.

\medskip

The symmetric group $\S_6$ acts on the 5-simplex by permuting its vertices, and
on its boundary which is the 4-sphere; composing the orientation-reversing
elements with ${\rm -id}_{\S^4}$ one obtains an orientation-preserving
orthogonal action of $\S_6$ on $S^4$.

\medskip

The finite subgroups of ${\rm SO}(4)$ and ${\rm O}(4)$ are well-known, see e.g.
[DV]. The group ${\rm SO}(4)$ is isomorphic to the central product $S^3
\times_{\Z_2} S^3$ of two unit quaternion groups $S^3$. There is a 2-fold
covering of Lie groups  $S^3 \to {\rm SO}(3)$ whose kernel is the central involution
of $S^3$. The finite subgroups of
${\rm SO}(3)$ are the polyhedral groups, that is cyclic $\Z_n$, dihedral $\D_{2n}$,
tetrahedral $\A_4$, octahedral $\S_4$ or dodecahedral $\A_5$; the finite subgroups
of $S^3$ are their preimages in $S^3$ which are cyclic, binary dihedral (or
generalized quaternion) $\D_{4n}^*$, binary tetrahedral $\A_4^*$,
binary octahedral $\S_4^*$ or binary dodecahedral $\A_5^*$.

\bigskip

A consequence of the Theorem and its proof is also the following characterization
of the finite subgroups of the orthogonal groups ${\rm SO}(5)$ and ${\rm
O}(5)$.

\bigskip

{\bf Corollary 2.} {\sl  Let $G$ be a finite subgroup of the orthogonal group
${\rm SO}(5)$  or  ${\rm O}(5)$.  Then one of the following cases occurs:

\medskip

i) $G$ is conjugate to  a subgroup of  $W = (\Z_2)^4 \rtimes \S_5$ or $\tilde W
= (\Z_2)^5 \rtimes \S_5$;

\medskip

 ii) $G$ is isomorphic to $\A_5$, $\S_5$, $\A_6$ or $\S_6$,   or
 to the product of one of these groups with  $\Z_2 = \langle {\rm
-id}_{\R^5} \rangle$;

\medskip

iii) $G$ is conjugate to a subgroup of  ${\rm O}(4)\times{\rm O}(1)$ or ${\rm
O}(3)\times{\rm O}(2)$.   }

\bigskip

Here ${\rm O}(4)\times{\rm O}(1)$ and ${\rm O}(3)\times{\rm O}(2)$ are
considered as subgroups of the orthogonal group ${\rm O}(5)$ in the natural
way. We note that $\A_5$ occurs as an irreducible subgroup of all three
orthogonal groups ${\rm SO}(3)$, ${\rm SO}(4)$ and ${\rm SO}(5)$, and in
particular the cases of Corollary 2 are again not mutually exclusive; see the
character tables in [C] for the irreducible representations of the groups in
ii).

\medskip

It should be noted that the proof of Corollary 2 is considerably easier than
the proof of the Theorem. For both the Theorem and Corollary 2 one has to
determine the finite simple groups which act on a homology 4-sphere resp. admit
an orthogonal action on the 4-sphere. In the case of the Theorem this is based
on [MeZ1, Theorem 1] and employs the Gorenstein-Harada classification of the
finite simple groups of sectional 2-rank at most four; for Corollary 2, this
can be replaced by much shorter arguments from the representation theory of
finite groups (see the Remark at the end of section 3).

\medskip

Finite groups acting on other types of 4-manifolds (including $S^2 \times S^2$
and  $\Bbb R^4$) are discussed in  [MeZ3] and [MeZ4],  see also the survey [E].

\bigskip

{\bf 2. Preliminaries}

\medskip

We introduce some notation and give some comments on the general structure of
the proof of the Theorem.  Let $G$ be a finite group acting on a homology
4-sphere as in the Theorem.  We shall consider the {\it maximal semisimple
normal subgroup} $E$ of $G$. Recall that a {\ semisimple group} is a central
product of quasisimple groups, and a {\it quasisimple group} a perfect central
extension of a simple group (see [S, chapter 6.6] or [KS, chapter 6.5]).
Crucial for the proof of the Theorem is [MeZ1, Lemma 4.3] which states that the
maximal semisimple normal group $E$ of $G$ is either trivial or isomorphic to
one of the following groups:
$$\A_5, \;\;\; \A_6, \;\;\; \A_5^*,  \;\;\; \A_5^* \times_{\Z_2} \A_5^*$$
($\A_5^* \times_{\Z_2} \A_5^*$ denotes the central product of two binary
dodecahedral groups $\A_5^*$); we note that this uses the Gorenstein-Harada
classification of the finite simple groups of sectional 2-rank at most four. In
the next section, we shall consider these different cases separately.  Note
that, if $E$ is nontrivial, we have to prove the stronger version iv') of the
Theorem.

\medskip

If $E$ is nontrivial we shall consider the  centralizer $C = C_G(E)$ of $E$ in
$G$ and the subgroup  $\tilde E$ of $G$ generated by $E$ and $C$; note that
$G/\tilde E$ is isomorphic to a subgroup of the outer automorphism group ${\rm
Out}(E) = {\rm Aut}(E)/{\rm Inn}(E)$ of $E$ which is small for the possible
groups $E$.

\medskip

If $E$ is trivial, the group $G$ may be solvable or nonsolvable. In this case,
we shall consider the {\it Fitting subgroup} $F$ of $G$ instead, that is the
maximal nilpotent normal subgroup of $G$. For an arbitrary finite group, the
subgroup $F^*$ generated  by  the Fitting subgroup $F$ and the maximal
semisimple normal subgroup $E$ is called the {\it generalized Fitting subgroup}
of $G$.  As its  main property, the generalized Fitting subgroup $F^*$ contains
its centralizer in $G$, hence  $G/F^*$ is isomorphic to a subgroup of the outer
automorphism group of $F^*$.  If $E$ is trivial, the generalized Fitting
subgroup $F^*$ coincides with the Fitting subgroup $F$; since $F$ is nilpotent
it is the direct product of its Sylow $p$-subgroups, for different primes $p$,
and each of these is normal in $G$.

\medskip

The following three Lemmas will be frequently used in the proof of the Theorem,
see [MeZ1, Lemma 4.1] for a proof of Lemma 1 (note that, in the statement of
case a) of Lemma 4.1 in [MeZ1], it is written erroneously 0-dimensional instead
of 2-dimensional). We note that, by general Smith fixed point theory, the fixed
point set of an orientation-preserving periodic diffeomorphism of prime power
order of a homology 4-sphere is either a 0-sphere or a 2-sphere (that is, a
homology sphere of even codimension).

\bigskip

{\bf Lemma 1.}  {\sl For a prime $p$, let $A$ be an elementary abelian
$p$-group acting orientation-preservingly on a homology 4-sphere.

a) If  $A$ has rank two and $p$ is odd, then $A$ contains exactly two cyclic
subgroups whose fixed point set is a 2-sphere, and the fixed point set of $A$
is a 0-sphere.

b) If $A$ has rank two and $p=2$, then $A$ contains at least one involution
whose fixed point set is a 2-sphere; in the case of three involutions with a
2-sphere as fixed point set the group $A$ has a 1-sphere as fixed point set.

c) If $A$ has rank three and $p=2$, then $A$ can contain either one or three
involutions whose fixed point set is  a 0-sphere.

d) If $A$ has rank four and $p=2$ , then $A$ contains exactly five involutions
whose fixed point set is  a 0-sphere.}

\bigskip

{\bf Lemma 2.}  {\sl Let $G$ be a finite group acting orientation-preservingly
on a homology 4-sphere. Suppose that  $G$ contains  a cyclic normal group $H$
of prime order $p$ such that the fixed point set of $H$ is a 0-sphere $S^0$;
then $G$ contains, of index at most two, a subgroup isomorphic to a subgroup of
${\rm SO}(4)$. Moreover, if $G$ acts orthogonally on $S^4$ then $G$ is
conjugate  to a subgroup of ${\rm O}(4)\times {\rm O}(1)$.}

\bigskip

{\it Proof.} The group $G$ leaves invariant the fixed point set $S^0$ of $H$
which consists of two points. A subgroup of index at most two fixes both points
and acts orthogonally and orientation-preservingly on a 3-sphere (the boundary
of a  regular invariant neighborhood of one of the two fixed points).

\medskip

If the action of $G$ is an orthogonal action on the 4-sphere then $G$
acts orthogonally on the equatorial 3-sphere of the 0-sphere $S^0$ and hence is a
subgroup of ${\rm O}(4)\times {\rm O}(1)$, up to conjugation.

\bigskip

{\bf Lemma 3.}  {\sl Let $G$ be a finite group acting orientation-preservingly
on a homology 4-sphere. Suppose that  $G$ contains  a cyclic normal group $H$
of prime order $p$ such that the fixed point set of $H$ is a 2-sphere $S^2$;
then $G$ is isomorphic to a subgroup of ${\rm O}(3) \times {\rm O}(2)$.
Moreover, if $G$ acts orthogonally on $S^4$ then $G$ is conjugate to a subgroup
of ${\rm O}(3) \times {\rm O}(2)$.}

\bigskip

{\it Proof.} The group $G$ leaves invariant the fixed point set $S^2$ of $H$. A
$G$-invariant regular neighbourhood of $S^2$ is diffeomorphic to the product of
$S^2$ with a 2-disk, so $G$ acts on its boundary $S^2 \times S^1$ (preserving its
Seifert fibration by circles). Now it is well-known that every finite group action
on $S^2 \times S^1$ is standard and conjugate to a subgroup of its isometry group
${\rm O}(3) \times {\rm O}(2)$.

\medskip

If $G$ acts orthogonally on $S^4$ then the group
$G$ leaves invariant $S^2$, the corresponding 3-dimensional subspace in $\Bbb{R}^5$
as well as its orthogonal complement, so up to conjugation  we are in ${\rm O}(3)
\times {\rm O}(2)$.

\bigskip

{\bf 3.  Proof of the Theorem  (and of Corollary 2)}

\medskip

We consider the various possibilities for the maximal normal semisimple
subgroup $E$ of $G$ as listed in section 2.

\bigskip

{\bf 3.1}  \hskip  2mm    Suppose that $E$ is isomorphic to $\A_5$.

\medskip

Let $C$ denote the  centralizer of $E$ in $G$ and $\tilde E$ the subgroup of
$G$ generated by $E$ and $C$. Then $G/\tilde E$ is isomorphic to a subgroup of
the outer automorphism group ${\rm Out}(E) = {\rm Aut}(E)/{\rm Inn}(E)$ of $E
\cong \A_5$ which has order two (see [C]).  If the centralizer $C$ of $E$ in
$G$ is trivial then $E =  \tilde E$ and $G$ is isomorphic to either  $\A_5$ or
$\S_5$ (case ii) of the Theorem).

\medskip

Suppose that $C$ is nontrivial.  We prove first that $C$ is cyclic or dihedral.
Let $S$ be a  Sylow 2-subgroup of $E$. The group $S$ is isomorphic to
$\Z_2\times \Z_2$ and its three involutions $t_1,t_2$ and $t_3$ are conjugate.
By Lemma 1, all three involutions have a 2-sphere as fixed point set and the
whole group $S$ has  a 1-sphere as fixed point set. Let $S^2$ denote the fixed
point set of $t_1$; note that $t_2$ and $t_3$ act as reflections in the same
1-sphere of $S^2$.

\medskip

Suppose that a nontrivial element $f$ in  $C$ acts trivially on $S^2$. Then $S$
and $f$ generate an abelian group and have a common fixed point on $S^2$; the
isotropy group of such a fixed point  has a faithful orthogonal representation
on euclidean space $\R^4$ (the tangent space of the fixed point). Since the
action of $G$ is orientation-preserving, the group generated by $S$ and $f$
acts faithfully on the linear subspace of dimension two orthogonal to the
tangent space of $S^ 2$; this gives a contradiction since only cyclic and
dihedral groups act orthogonally and faithfully on $\R^2$.

\medskip

So it follows that the action of $C$ on $S^2$ is faithful. We can also suppose
that the action of $C$ on $S^2$ is orientation-preserving, otherwise we compose
the orientation-reversing elements of $C$ with $t_2$ obtaining a new group
isomorphic to $C$ and with an orientation-preserving action. We consider the
finite groups which admit a faithful,  orientation-preserving action on the
2-sphere; $C$ cannot be isomorphic to $\A_4$, $\S_4$ or $\A_5$ since the
reflections $t_2$ and $t_3$ of $S^2$ in a circle  do not commute with such a
group. So we conclude that $C$ is cyclic or dihedral.

\medskip

Suppose that $C$ contains  a cyclic subgroup $C'$ of prime order $p$ which is
normal in $G$; this is always the case if  $C$ is either cyclic, or dihedral of
order strictly greater than four. By Smith fixed point theory,  the fixed point
set of $C'$ is either a 0-sphere or a 2-sphere.

\medskip

If the fixed point set of $C'$ is a 2-sphere $S^2$, we shall prove that we are
in  case iii) of the Theorem. We consider the  subgroup $C''$ consisting of
elements of $G$ which act trivially on $S^2$; it is normal and cyclic (it acts
as a group of rotations around $S^ 2$). The group $E\cong\A_5$ acts by
conjugation on $C''$; since $E$ is simple and the automorphism group of $C''$
is abelian, $E$ acts trivially on $C''$. If there is an element $g$ in $G$ not
contained in $\tilde{E}$, then $g$ acts nontrivially on $E$ and in particular
$g$ is not contained in  $C''$; this implies that the factor group $G/C''$,
which acts faithfully on $S^2$, contains a subgroup isomorphic to the symmetric
group $\S_5$. However $\S_5$ does not act faithfully on the 2-sphere, so it
follows that $G$ and $\tilde E$ coincide, $G$ is isomorphic to $E\times C''$
and we are in case iii) of the Theorem.

\medskip

On the other hand, if the fixed point set of $C'$ is a 0-sphere $S^0$ that is
consists of two points, we shall prove that we are in case iv) of  the Theorem.
The group $G$ contains $G_0$, a subgroup  of index at most two that fixes both
points; this subgroup acts faithfully and orientation-preservingly on a
3-sphere which is the boundary of a  regular invariant neighborhood of one of
the two fixed points. By the classification of finite subgroups of ${\rm
SO}(4)$  in [DV], it follows that $G_0$ is isomorphic to $\A_5\times \Z_2$; the
fixed point set of the involution in the center of $G_0$ is $S^0$. If $G_0=G$
we are done. Suppose that  $G_0$ has index two, let $f$ be an element in $G$
but not in $G_0$. Since the automorphism group of $\A_5$ is isomorphic to
$\S_5$ (see [C]), we can suppose that $f^2$ acts trivially on $E\cong \A_5$ and
$f^2$ is in the center of $G_0$. By the Lefschetz fixed point theorem $f$ has
non-empty fixed point set; if $f^2$  is the involution in the center of $G_0$,
the fixed point set of $f$ is $S^0$ and this is impossible. We can suppose that
$f$ is an involution. If $f$ does not act trivially on $E$, we obtain that
$G\cong \S_5 \times \Z_2$ that is a subgroup of ${\rm O}(4)$  (case 51 in
[DV]). If $f$ acts trivially on $E$ we have $G\cong\A_5\times \Z_2\times \Z_2$
and this is again a subgroup of ${\rm O}(4)$ (case 49 in [DV]).

\medskip

Finally we can  suppose that $C\cong \Z_2\times \Z_2$ is dihedral of order
four. By Lemma 1 we have three possibilities for the fixed point sets of the
involutions in $C$: the fixed point set of each involution is a 2-spheres, or
the fixed point set of exactly one involution is a 0-sphere, or the fixed point
set of exactly one involution is a 2-sphere. In the last two possibilities  an
involution in $C$ is central in $G$ (it is not conjugate to any other element
in $C$) and we are   in the previous case. If the fixed point set of each
involution is a 2-sphere, by Lemma 1 the global fixed point set of $C$ is a
1-sphere $S^1$. The subgroup  consisting of elements of $G$ acting trivially on
$S^1$ is normal in $G$. Since $\A_5$ does not act faithfully on a 1-sphere, the
action of $E$ and hence also of $\tilde E$ on $S^1$ has  to be trivial. Then
$\tilde E = E \times C$ is contained in the isotropy group of any point of
$S^1$ and acts faithfully and orthogonally on $\R^4$ (the tangent space of the
point). The group  $E \times C$ fixes pointwise the tangent space  of $S^1$
and hence  acts faithfully and orthogonally also on the orthogonal subspace  of
dimension 3. But $E \times C \cong \A_5 \times \Z_2 \times \Z_2$ does not act
faithfully and orthogonally on $\R^3$ (or $S^2$). This contradiction completes
the proof in case  3.1.

\bigskip

{\bf 3.2}  \hskip 2mm  Suppose next that  $E$ is isomorphic to  either $\A_5^*$  or
$\A_5^*  \times_{\Z_2} \A_5^*$.

\medskip

By Smith  fixed point theory, the fixed point set of
the central involution $t$ of $E$ is either a 0-sphere or a 2-sphere.

\medskip

Suppose that the fixed point set of $t$ is a 2-sphere $S^2$; we shall show that
this case really does not occur. Note that $E$ cannot be isomorphic to $\A_5^*
\times_{\Z_2} \A_5^*$ since no quotient group of $\A_5^* \times_{\Z_2} \A_5^*$
by a cyclic subgroup (the group fixing pointwise $S^2$) admits an action on a
2-sphere.  So there remains the case  $E \cong \A_5^*$. Then $E$ admits a
faithful action on a regular neighborhood of $S^2$ which is diffeomorphic to the
product of $S^2$ with a 2-disk (since the self-intersection
number of $S^2$ in a homology 4-sphere has to be zero), and on its boundary
$S^2  \times S^1$ (preserving its Seifert fibration by circles). However, it is
well-known that every action of a finite group on $S^2 \times S^1$ is standard
(preserves the product structure, up to conjugation), and in particular that
$\A_5^*$ does not admit a faithful action on $S^2 \times S^1$. So also this
case does not occur.

\medskip

On the other hand, if the fixed point set of $t$  is a 0-sphere $S^0$ we shall
prove that we are in case iv) of the Theorem.  A subgroup $G_0$ of index two of
$G$ fixes both points of $S^0$ and is isomorphic to a subgroup of  ${\rm
SO}(4)$.   If $G=G_0$ we are done, so we can suppose  that $G_0$ has index two
in $G$. Let $f$ be an element in $G$ but not in  $G_0$.

\medskip

Suppose first that $f$ normalizes  $G_1$, a   subgroup of $G_0$ isomorphic to
$\A^*_5$  containing $t$ in its center;  this is always the case if
$E\cong\A^*_5$. The group $G_1$ fixes both points and $G_1$  acts orthogonally
and orientation-preservingly on a 3-sphere which is the boundary of a  regular
invariant neighborhood of one of the two fixed points set. The action of
$\A^*_5$ on $S^3$ has to be free, so the fixed point set of any element in
$G_1$ is $S^0$. Up to composition by an element of $G_1$, we can suppose either
that $f$  centralizes $G_1$ or that $f$ is represented in ${\rm Aut}(G_1)\cong
\S_5$ by a transposition; in any case $f$ centralizes $h$, an element of order
three in $G_1$. We can suppose that the order of $f$ is a power $2^m$ of two.
Since $f \notin G_0$, $f$ exchanges the two points of $S^0$, the fixed point
set of $h$. By the Lefschetz fixed point theorem $fh$ has a fixed point not in
$S^0$; then the same holds for $(fh)^{2^m} = f^{2^m}h^{2^m} = h^{\pm 1}$ which
is a contradiction, so this case cannot occur.

\medskip

We can suppose that $E$ is isomorphic to   $\A_5^*\times_{\Z_2} \A_5^*$ and $f$
exchanges by conjugation the two quasisimple components. Since
$\A_5^*\times_{\Z_2} \A_5^*$ is maximal among the groups acting
orientation-preservingly on a 3-sphere, the subgroup $G_0$ coincides with $E$.
We describe  $E\cong \A_5^*\times_{\Z_2} \A_5^*$ as the factor group of
$\A^*_5\times \A^*_5=\{(x,y)| x,y\in\A^*_5\}$ by the normal subgroup
$\{(c,c),(1,1)\}$ where $c$ is the central involution in  $\A^*_5$. By [S,
Theorem 6.11]  the automorphism group of $E$ is isomorphic to a semidirect
product of the normal subgroup   ${\rm Aut}(\A_5)\times {\rm Aut}(\A_5)\cong
\S_5 \times \S_5$ with a group of order two generated by the automorphism
$\phi(x,y)=(y,x)$.  Since $f^2\in E$ and $f$ exchanges the two components, up
to composition of $f$ by an element of $E$, we can suppose that the
automorphism induced by $f$ is either $\phi(x,y)=(y,x)$ or
$\phi'(x,y)=(\sigma(y),\sigma(x))$ where $\sigma$ is a non-inner automorphism
of order two of $\A^*_5$. In any  case we can suppose that $f^2$ acts trivially
by conjugation on  $E$, that is $f^2$ is contained in the center of $E$. By the
Lefschetz fixed point theorem $f$ has non empty fixed point set. If $f^2$ is
the non trivial element in the center of $E$,  the fixed point set of $f$
coincides with $S^0$, the fixed point set of the involution in the center of
$E$, and this is impossible.  We obtain that $f$ is an involution and that $G$
is a splitting extension of $E$ by the subgroup of order two generated by $f$.
The two possible automorphisms induced by $f$, that are $\phi$ and $\phi'$, are
conjugate in the automorphism group of $E$, hence the possible extension of $E$
by $f$ is unique up to isomorphism. In this case $G$ has to be  isomorphic to
the unique extension of $\A_5^*\times_{\Z_2} \A_5^*$ that appears in the list
of finite subgroup of  ${\rm O}(4)$ (the case  50 in [DV]).

\bigskip

{\bf 3.3}  \hskip 2mm  The case $E \cong \A_6$ is considered in [MeZ1, Theorem 2.1],
and we will not repeat the arguments here; one shows that $G$ is isomorphic to
either  $\A_6$ or $\S_6$, so we are in case ii) of  the Theorem.

\medskip

Finally we come to the last case:

\bigskip

{\bf 3.4}  \hskip 2mm   Suppose that $E$ is trivial.

\medskip

Now $G$ may be solvable or nonsolvable; if $G$ is nonsolvable we
have to prove the stronger version iv') of the Theorem. In the following, we
start by proving version iv) of the Theorem, and comment on the nonsolvable case in
the appropriate instance.

\medskip

Let $F$ denote the Fitting subgroup of $G$; since $E$ is
trivial, $F$ coincides with the generalized Fitting subgroup and contains
its centralizer in $G$ (see section 2).  Also, the nilpotent group $F$ is the direct
product of its Sylow  $p$-subgroups, for different primes $p$, and each of these is
normal in $G$.

\medskip

Suppose that  $F$ of $G$ contains a nontrivial p-Sylow subgroup
$S_p$ with $p$ odd. We consider the maximal elementary abelian subgroup $Z$ in the
center of $S_p$; the subgroup $Z$ is normal in $G$. By Smith fixed point theory
(see [MeZ1, Lemma 2.3]), the rank of $Z$ is at most two. If the rank is
one, either Lemma 2 or Lemma 3 applies and we are done. If the rank of $Z$
is two, by Lemma 1 the fixed point set of $Z$ is a 0-sphere $S^0$. The
group $G$ leaves invariant the 0-sphere and we can conclude as in Lemma 2.

\medskip

We can suppose hence that the Fitting subgroup $F$ is a $2$-group. Let $Z$ denote
again the maximal elementary abelian subgroup contained in the center of $F$;
the subgroup $Z_1$ is normal in $G$.  By Smith theory (see [MeZ1, Lemma 2.3.]),
the rank of $Z_1$ is at most four, and accordingly we consider the following
four cases.

\medskip

If $Z$ is cyclic we are done by Lemma 2 or 3.

\medskip

Suppose that $Z$ has rank two. The group  $Z$ contains three involutions. If
exactly one involution in $Z$ has a 2-sphere or a 0-sphere as fixed point set,
$G$ has again a cyclic normal subgroup and Lemma 2 or Lemma 3 applies. By Lemma
1 the unique other possibility is that all three involutions have a 2-sphere as
fixed point set; in this case the fixed point set of $Z$  is a 1-sphere $S^1$,
invariant under the actions of $G$.  The boundary of a $G$-invariant regular
neighbourhood of $S^1$ is diffeomorphic to $S^2 \times S^1$. Since finite group
actions on $S^2 \times S^1$ are standard, $G$ is isomorphic to a subgroup of
its isometry group ${\rm O}(3)\times {\rm O}(2)$.   Moreover, if $G$ acts
orthogonally on $S^4$ then it leaves invariant $S^1$, the corresponding
subspace of dimension two in $\Bbb{R}^5$ and its orthogonal complement, hence
is conjugate to a finite subgroup of ${\rm O}(3)\times {\rm O}(2)$.

\medskip

If $Z$ has rank three then by Lemma 1 it has either one or three involutions
with a 0-sphere  as fixed point set. In the first case we have a cyclic normal
subgroup and Lemma 2  applies. In the second case the three involutions with a
0-sphere as fixed point set generate $Z$ since, by Lemma 1, a subgroup of
rank two cannot contain three involutions with a 0-sphere as fixed point set.
The group $G$  permutes by conjugation these three involutions, and their
product is an involution which is central in $G$ and has a 2-sphere as fixed
point set; in this case we are done by Lemma 3.

\medskip

Finally we suppose that $Z$ has rank four. We note that, by [MeZ1, Lemma 4.2],
we are exactly in this situation if $G$ is nonsolvable (that is, in all cases
of 3.4 considered so far the group $G$ is solvable and version iv) of the
Theorem applies).

\medskip

By Lemma 1  the group  $Z$ contains at least one involution $t$ which
has a 0-sphere
$S^0$ as fixed point set. The group $F$ leaves invariant $S^0$ and a subgroup $F_0$
of index at most two of $F$ fixes both points of $S^0$. The group $F_0$ acts
faithfully on a 3-sphere that is the boundary of a regular invariant neighborhood
of one of the two fixed points. Moreover $F_0$ contains in its center an elementary
abelian 2-subgroup of rank three (the group $F_0 \cap Z$); since the only 2-group
acting on the 3-sphere with this property is $(\Bbb{Z}_2)^3$ (see [MeZ2,
Propositions 2 and 3]), it follows that the Fitting
subgroup $F$ coincides with $Z$ and is an elementary abelian 2-group of
rank four.

\medskip

Recall that $F$ contains its centralizer in $G$. By Lemma 1, $F$ contains
exactly five involutions with 0-dimensional fixed point set; they generate the
group $F$ because by Lemma 1 the  subgroups of index two  contain at most three
such involutions. We conclude that $G/F$ acts faithfully on the set of the five
involutions of $F$ with 0-dimensional fixed point set,  and hence $G/F$ is
isomorphic to a subgroup of $\Bbb{S}_5$.

\medskip

Suppose first that $G$ is solvable. Considering the solvable subgroups of
$\Bbb{S}_5$, one can distinguish two cases then.  If there is a subgroup $\Z_2$ or
$(\Z_2)^2$ of $Z \cong (\Z_2)^4$, invariant under the action of $G/F$ and hence
normal in $G$, then one concludes as above that  cases iii) or iv) of the Theorem
apply. The only other possibility is that $G/F$ is cyclic of order five or
dihedral of order ten, and then it is easy to see that $G$ is the semidirect product
of $Z$ with $G/F$ and hence isomorphic to a subgroup of $W = (\Z_2)^4
\rtimes \S_5$.

\medskip

On the other hand, if $G$ is nonsolvable, then
$G/F$ is isomorphic to $\S_5$ or $\A_5$. Moreover, the action of $\S_5$ or $\A_5$
on the five involutions with 0-dimensional fixed point set is the standard
permutation action, and the action of $\S_5$ or $\A_5$ on $F \cong (\Z_2)^4$
coincides with the action in the semidirect products  $W =
(\Z_2)^4 \rtimes \S_5$ or its subgroup $W_0 = (\Z_2)^4 \rtimes \A_5$ of index two.
So it remains to show that $G$ is in fact such a semidirect product (that is, a
split extension).

\medskip

Suppose first that $G/F$ is isomorphic to $\A_5$. Then $G$ is a perfect group
of order 960. The perfect groups of small order have been determined in [HP]
(implemented also in the group theory software GAP), and there are exactly two
perfect groups of order 960 ([HP, cases 4.1 and 4.2, p.118]). Both are
semidirect products $(\Z_2)^4 \rtimes \A_5$, one is isomorphic to the subgroup
$W_0$ of index two of $W$, so $G$ is isomorphic to $W_0$ in
this case.  (We note that there are exactly two subgroups $\A_5$ of the automorphism
group ${\rm GL}(4,2) \cong \A_8$ of $(\Z_2)^4$, up to conjugation, corresponding to
the two actions of $\A_5$ on $(\Z_2)^4$, see also [HP, case 0.1, p.116]. One may
then consider the second cohomology $H^2(\A_5; (\Z_2)^4)$ which is trivial for both
actions ([HP, p.118]); since $H^2(\A_5; (\Z_2)^4)$ classifies extension of
$(\Z_2)^4$ with factor group $\A_5$ ([Bro; Theorem IV.3.12]), the only possible
extensions are the two semidirect products.)

\medskip

Now suppose that $G/F$ is isomorphic to $\S_5$. By the first case, $G$ contains
$W_0 = (\Z_2)^4 \rtimes \A_5$ as a subgroup of index two. Now there is exactly
one embedding of $\A_5$ into $W_0$, up to conjugation: this can be checked
directly (e.g. by GAP), or it follows from the fact that  $H^1(\A_5; (\Z_2)^4)$
is trivial for the action corresponding to $W_0$ ([HP, p.118]) (the first
cohomology $H^1(\A_5; (\Z_2)^4)$ classifies injections of $\A_5$ into the
semidirect product up to conjugation, see [Bro, Proposition IV.2.3]).  Let $x$
be a coset representative of $W_0$ in $G$, representing an element of order two
in the factor group $\S_5$ of $G$, so $x^2 \in (\Z_2)^4$.  By the preceding, we
can assume that $x \A_5 x^{-1} = \A_5$. Then necessarily $x^2 = 1$ (since
otherwise $x^2 \A_5 x^{-2} \ne \A_5$), and hence $\A_5$ and $x$ generate a
subgroup $\S_5$ of $G$ and $G$ splits as a semidirect product $(\Z_2)^4 \rtimes
\S_5$ isomorphic to $W$.

\bigskip

This completes the proof of the Theorem.

\bigskip

Concerning the proof of Corollary 2 in the last case  $F \cong (\Z_2)^4$,
suppose that  $G$ acts orthogonally on $S^4$  and  consider the corresponding
orthogonal action of $G$ on $\Bbb{R}^5$. In this case the five  1-dimensional
subspaces of
$\Bbb{R}^5$ which are the  fixed point sets of the five involutions  in $F$ are
pairwise orthogonal and, up to conjugation, we can suppose that $F$ is the
orientation-preserving subgroup of the group generated by the inversion of the
coordinates in $\Bbb{R}^5$ (every elementary abelian 2-subgroup of an
orthogonal group is diagonalizable).  Any element in the normalizer of $F$ in
${\rm O}(5)$ leaves invariant the set of five 1-dimensional subspaces which are
the fixed point sets of the five involutions in $F$, and the elements leaving
invariant each of these 1-dimensional subspaces are contained in $F$. The
normalizer  of $F$ in ${\rm O}(5)$  is generated by $F$ and by the group of
permutations of the coordinates, and in particular  $G$ is conjugated to a
subgroup of the group $W$.

\bigskip

As noted in the introduction, for the proof of Corollary 2 the
Gorenstein-Harada classification of the finite simple groups of sectional
2-rank at most four can be replaced by arguments from the representation theory
of finite groups. In fact, the finite groups which admit a faithful,
irreducible representation in degree five (or equivalently, are irreducible
subgroups of the linear group ${\rm SL}(5,\Bbb C)$) have been determined in
[Bra], and the simple groups occurring are the alternating groups $\A_5$ and
$\A_6$, the linear fractional group ${\rm PSL}(2,11)$ and the unitary or
symplectic group ${\rm PSU}(4,2) \cong {\rm PSp}(4,3)$ (or, in another
notation, the group ${\rm U}_4(2) \cong {\rm S}_4(3)$); in dimensions less than
five there occurs in addition the linear fractional group ${\rm PSL}(2,7)$,
with an irreducible representation in dimension three. Since none of these
groups except $\A_5$ and $\A_6$ admits a faithful {\it real} representation in
dimension five (see [C]), we are left with the groups $\A_5$ and $\A_6$.

\bigskip

{\bf 4.  Some comments on the situation in dimension three}

\medskip

We close with a short discussion of the class of finite, in particular finite
nonsolvable groups $G$ which admit a faithful, smooth, orientation-preserving
action on a homology 3-sphere.

\medskip

We consider the class of groups $Q(8a,b,c)$ in [Mn] which have periodic
cohomology of period four but do not admit a faithful, linear action on $S^3$.
We will assume in the following that $a > b > c \ge 1$ are odd coprime
integers; then $Q(8a,b,c)$ is a semidirect product $\Z_{abc} \rtimes Q_8$ of a
normal cyclic subgroup  $\Z_a \times \Z_b \times \Z_c  \cong \Z_{abc}$ by the
quaternion group $Q_8 = \{\pm 1,\pm i,\pm j,\pm k\} \cong \D_8^*$ of order
eight, where $i,j$ and $k$ act trivially on $\Z_a, \Z_b$ and $\Z_c$,
respectively, and in a dihedral way on the other two.

\medskip

It has been shown by Milgram [Mg] that some of the groups $Q(8a,b,c)$ admit a
faithful, free (in particular orientation-preserving) action on a homology
3-sphere (and some others do not; see also the comments in [K, p.173, Update A
to Problem 3.37]). On the other hand, it is not difficult to see that none of
the groups $Q(8a,b,c)$ admits a faithful, linear action on the 3-sphere, i.e.
they are not isomorphic to subgroups of the orthogonal group ${\rm O}(4)$.

\medskip

The groups $Q(8a,b,c)$ are clearly solvable; concerning nonsolvable groups, it
is shown in [Z] that the finite nonsolvable groups which admit a faithful,
smooth, orientation-preserving action on a homology 3-sphere are exactly the
finite nonsolvable subgroups of the orthogonal group ${\rm SO}(4)  \cong S^3
\times _{\Z_2} S^3$, plus possibly two other types of groups which are the
central products
$$\A_5^* \times _{\Z_2}Q(8a,b,c)$$ and their subgroups
$$\A_5^* \times _{\Z_2} \D_{4a}^* \times \Z_b.$$ In turn these have a subgroup
$$\D_8^* \times _{\Z_2} \D_{4a}^* \times \Z_b,$$ and it is easy to see that, for
odd, coprime integers $a,b \ge 3$, a group $\D_8^* \times _{\Z_2} \D_{4a}^*
\times \Z_b$ does not admit a faithful, linear, orientation-preserving action
on $S^3$. So this leads naturally to the following

\bigskip

{\bf Question.}  Suppose $a$ and $b$ are odd coprime integers grater then one.
Does
$$\D_8^* \times_{\Z_2} \D_{4a}^* \times \Z_b$$ admit a faithful,
orientation-preserving action on a homology 3-sphere? (If $a$ is even then
there no such action exists by [Z, Lemma].)

\medskip

If the answer is no (as we expect) then by [Z] the class of nonsolvable groups
admitting a faithful, smooth, orientation-preserving action on a homology
3-sphere coincides exactly with the class of nonsolvable subgroups of the
orthogonal SO(4); otherwise, one has some new and easy solvable groups which
admit a faithful, orientation-preserving action on a homology 3-sphere but are
not subgroups of SO(4) (similar as for the Milnor groups $Q(8a,b,c)$; but
whereas a Milnor group can admit only {\it free} faithful actions on a homology
3-sphere, a faithful action of a group $\D_8^* \times_{\Z_2} \D_{4a}^* \times
\Z_b$ is necessarily nonfree (since it has a subgroup $\Z_2 \times \Z_2$)).

\bigskip \bigskip

\centerline {\bf References}

\bigskip

\item {[Bra]}  R. Brauer, {\it \"Uber endliche lineare Gruppen von Primzahlgrad.}
Math. Ann. 169, 73-96  (1967)

\item {[Bre]} G. Bredon, {\it Introduction to Compact Transformation Groups.}
Academic Press, New York 1972

\item {[Bro]} K.S. Brown, {\it Cohomology of Groups.}  Graduate Texts in Mathematics
87, Springer 1982

\item {[C]} J.H. Conway, R.T. Curtis, S.P. Norton, R.A. Parker, R.A. Wilson, {\it Atlas
of Finite Groups.} Oxford University Press 1985

\item {[DV]}  P. Du Val, {\it  Homographies, Quaternions and Rotations.} Oxford
Math. Monographs, Oxford University Press 1964

\item {[E]} A.L. Edmonds, {\it A survey of group actions on 4-manifolds.}
arXiv:0907.0454

\item {[HP]}  D.F. Holt,  W. Plesken, {\it Perfect Groups.} Oxford University Press
1989

\item {[K]} R. Kirby, {\it Problems in low-dimensional topology.}  Geometric
Topology. AMS/IP Studies in Advanced Mathematics Volume 2, part 2, 35-358
(1997)

\item {[KS]} H. Kurzweil, B. Stellmacher, {\it The Theory of Finite Groups.}
Universitext, Springer 2004

\item {[MeZ1]} M. Mecchia, B. Zimmermann, {\it On finite simple and nonsolvable
groups acting on homology 4-spheres.} Top. Appl. 153,  2933-2942  (2006)

\item {[MeZ2]} M. Mecchia, B. Zimmermann, {\it On finite groups acting on
$\Z_2$-homology 3-spheres.} Math. Z. 248, 675-693 (2004)

\item {[MeZ3]} M. Mecchia, B. Zimmermann, {\it On finite simple and nonsolvable groups
acting on closed 4-manifolds.}  Pac. J. Math. 243, 357-374 (2009)
(arXiv:0803.4454)

\item {[MeZ4]} M. Mecchia, B. Zimmermann, {\it On finite groups acting on acyclic
low-dimensional manifolds.}  arXiv:0808.0999v2

\item {[Mg]} R.J. Milgram, {\it Evaluating the Swan finiteness obstruction for
finite groups.} Algebraic and Geometric Topology. Lecture Notes in Math. 1126
(Springer 1985), 127-158

\item {[Mn]} J. Milnor, {\it Groups which act on $S^n$ without fixed points.} Amer.
J. Math. 79, 623-630  (1957)

\item {[S]} M. Suzuki, {\it Group Theory II.}  Springer-Verlag 1982

\item {[Z]} B. Zimmermann, {\it On the classification of finite groups acting on
homology 3-spheres.} Pacific J. Math. 217, 387-395 (2004)

\bye